\documentclass[a4paper,reqno,11pt]{amsart}

\usepackage{amsmath}
\usepackage{amssymb}
\usepackage{amsthm}
\usepackage{mathtools}
\usepackage{calc}
\usepackage{enumitem}
\usepackage{fullpage}
\usepackage[nocompress]{cite}
\usepackage{thmtools, thm-restate}
\usepackage[hidelinks]{hyperref}
\usepackage{cleveref}
\usepackage{graphicx}

\newtheorem{theorem}{Theorem}[section]

\newtheorem{lemma}[theorem]{Lemma}

\newtheorem{observation}[theorem]{Observation}

\newtheorem{claim}[theorem]{Claim}
\newtheorem*{claim*}{Claim}

\theoremstyle{definition}

\newtheorem{question}[theorem]{Question}
\newtheorem*{qu*}{Question}
\theoremstyle{remark}

\newcommand{\nlongrightarrow}{%
  \mathrel{\ooalign{%
    $\longrightarrow$\cr
    \hfil\scalebox{0.6}{$/$}\hfil\cr
  }}%
}

\newcommand{\nindarrow}{\overset{\mathrm{ind}}{\nlongrightarrow}}

\newcommand{\indarrow}{\overset{\mathrm{ind}}{\longrightarrow}}

\newcommand\N{\mathbb{N}}

\newcommand\cA{\mathcal{A}}

\newcommand\cC{\mathcal{C}}

\newcommand\cH{\mathcal{H}}
\newcommand\cI{\mathcal{I}}

\newcommand\cK{\mathcal{K}}
\newcommand\cL{\mathcal{L}}

\newcommand\cS{\mathcal{S}}

\renewcommand\Pr{\operatorname{\mathbb{P}}}

\newcommand\eps{\varepsilon}
\renewcommand\leq{\leqslant}
\renewcommand\geq{\geqslant}
\renewcommand\le{\leqslant}
\renewcommand\ge{\geqslant}
\renewcommand\to{\rightarrow}

\newcommand\F{\mathbb F}
\newcommand\Fq{\mathbb F_q}
\newcommand\Ex{\mathbb{E}}
\newcommand\Var{\mathrm{Var}}

\crefname{thm}{theorem}{theorems}
\crefname{lem}{lemma}{lemmas}
\crefname{obs}{observation}{observations}
\crefname{defi}{definition}{definitions}
\Crefname{thm}{Theorem}{Theorems}
\Crefname{prop}{Proposition}{Propositions}
\Crefname{obs}{Observation}{Observations}
\Crefname{claim}{Claim}{Claims}
\Crefname{defi}{Definition}{Definitions}
\crefname{subsection}{subsection}{subsections}

\crefformat{enumi}{#2item~#1#3}
\Crefformat{enumi}{#2Item~#1#3}
\crefmultiformat{enumi}
  {items~#2#1#3}
  { and~#2#1#3}
  {, #2#1#3}
  {, and~#2#1#3}
\crefrangeformat{enumi}
  {items~#3#1#4--#5#2#6}

\title{Off-Diagonal Induced Ramsey Numbers}
\author{Bruno Andrades, Marcelo Campos, Rafael Filipe, Anastasiia Kats}

\address{IMPA, Estrada Dona Castorina 110, Jardim Bot\^anico, Rio de Janeiro, 22460-320, Brasil}

\email{bruno.andrades|marcelo.campos|rafael.santos|kats.anastasiia@impa.br}

\thanks{BA and RF were supported by CNPq; MC was partially supported by Serrapilheira (grant R-2412-51283); and AK was supported by a Bolsa do IMPA}

\begin{document}

\begin{abstract}
Given graphs $F$ and $H$, the induced Ramsey number $r_{\mathrm{ind}}(F,H)$ is the smallest integer $n$ for which there exists a graph $G$ on $n$ vertices such that every red-blue colouring of $E(G)$ contains either a blue induced copy of $F$ or a red induced copy of $H$. A conjecture of Kohayakawa, Pr\"omel and R\"odl asserts that, for every fixed graph $F$, this number is polynomial in $v(H)$.
We prove this conjecture, showing that there is an absolute constant $C > 0$ such that
$$r_{\mathrm{ind}}(F,H) \leq v(H)^{Cv(F)^2}$$
for all graphs $F$ and $H$. The proof combines a variant of a random geometric construction of Kohayakawa, Pr\"omel and R\"odl with the hypergraph container method.
\end{abstract}

\maketitle

\section{Introduction}

The induced Ramsey number $r_{\mathrm{ind}}(H)$ of a graph $H$ is the minimum number of vertices of a graph $G$ with the property that in every red-blue colouring of the edges of $G$, there is an induced copy of $H$ all of whose edges are the same colour. These numbers were shown to be finite for all $H$ by Deuber~\cite{Deuber75}, Erd\H{o}s, Hajnal and P\'osa~\cite{EHP75}, and R\"odl~\cite{Rodl73}, but the bounds given by these initial proofs were only double exponential (or worse) in the number of vertices of $H$. Following significant improvements by Kohayakawa, Pr\"omel and R\"odl~\cite{KPR98} in the 1990s, and by Conlon, Fox and Sudakov~\cite{CFS12} in the 2010s, an exponential upper bound 
$$r_{\mathrm{ind}}(H) \le 2^{O(k)}$$
for every graph $H$ with $k$ vertices was finally shown by Arag\~ao, Campos, Dahia, Filipe and Marciano~\cite{ACDFM25}, confirming a conjecture of Erd\H{o}s~\cite{Erdos84}.

In this paper we will study the following natural asymmetric version of this question, which was first considered by Kohayakawa, Pr\"omel and R\"odl~\cite{KPR98}. Given two graphs $F$ and $H$, we write 
$$G \indarrow (F,H)$$
if for every red-blue colouring of the edges of $G$, there exists either a blue induced copy of $F$ (that is, an induced copy of $F$ in $G$ all of whose edges are blue), or a red induced copy of $H$. The \emph{induced Ramsey number} of the pair ($F,H)$ is
$$r_{\mathrm{ind}}(F,H) = \min\big\{ v(G) : G \indarrow (F,H) \big\}.$$
In particular, note that $r_\mathrm{ind}(K_s,K_t) = R(s,t)$, the usual off-diagonal Ramsey number, since every clique in a graph $G$ is also an induced subgraph of $G$.

Kohayakawa, Pr\"omel and R\"odl~\cite{KPR98} proved that 
\begin{equation}\label{eq:KPR}
r_{\mathrm{ind}}(F,H) \le v(H)^{O(v(F) \log \chi(H))}
\end{equation}
for all graphs $F$ and $H$ with at least one edge, and moreover conjectured~\cite[Conjecture~2]{KPR98} that $r_{\mathrm{ind}}(F,H)$ is bounded by a polynomial function of $v(H)$, for every fixed graph $F$. That is, for each graph $F$, there exists $C(F) > 0$ such that
$$r_{\mathrm{ind}}(F,H) \leq v(H)^{C(F)}$$
for every graph $H$. Our main result confirms this conjecture with $C(F) = O(v(F)^2)$. 

\begin{theorem}\label{thm:main}
There exists an absolute constant $C > 0$ such that  
\begin{equation}\label{eq:rFH:bound}
r_{\mathrm{ind}}(F,H) \le v(H)^{C v(F)^2}
\end{equation}
for all graphs $F$ and $H$. 
\end{theorem}

We also prove the following multicolour version of \Cref{thm:main}. Given graphs
$F_1,\ldots,F_r$, let $r_{\mathrm{ind}}(F_1,\ldots,F_r)$ denote the smallest number of vertices of a graph $G$ with the property that every $r$-colouring of $E(G)$ contains an induced copy of $F_i$ in colour $i$ for some $i \in [r]$.

\begin{theorem}\label{thm:multicolour}
There exists an absolute constant $C>0$ such that for every $r \in \N$
$$r_{\mathrm{ind}}(F_1,\ldots,F_r,H) \le v(H)^{C r^2 s^2}$$
for all graphs $F_1,\ldots,F_r$ and $H$ with $v(F_1) = \cdots = v(F_r) = s$ and $v(H) \ge Cr$. 
\end{theorem}

The construction we will use to prove~\Cref{thm:main,thm:multicolour} is based on that of Kohayakawa, Pr\"omel and R\"odl~\cite{KPR98}, who proved the bound~\eqref{eq:KPR} by placing a random blow-up of $H$ in each line of a projective plane $P$. Note that this graph has many induced copies of $H$, since each pair of vertices is in a unique line of $P$, so the copies do not interfere with one another. It is also not hard to show (see \Cref{obs:ksquared}) that in any red-blue colouring of a balanced blow-up of $H$ with no induced red copy of $H$, the density of blue edges is at least $v(H)^{-2}$. In~\cite{KPR98}, the authors used this fact to show that the blue graph is pseudorandom, and hence to embed $F$ in it. 

Our strategy will have some high level similarities, but also two important changes. First, our host graph $G$ will be slightly different: instead of using a projective plane, we will place a random blow-up of $H$ in each line of $\F_q^d$, where $d = 4s$. Working in a higher dimensional space will guarantee more randomness in our construction, which will be important in the proof.

The second change is that, instead of using a pseudorandom property to embed $F$ into the blue graph, we will exploit the randomness of the host graph and use the method of hypergraph containers. This method, developed independently by Balogh, Morris, and Samotij~\cite{BMS15} and by Saxton and Thomason~\cite{ST15}, has become a powerful and versatile tool for studying families of sets that avoid prescribed configurations; see~\cite{S15, BMS18survey} for surveys. 

The container method was first applied to induced Ramsey numbers by Conlon, Dellamonica, La Fleur, Rödl and Schacht~\cite{CDFRS17}, who used it to prove an upper bound for multicolour induced Ramsey numbers of uniform hypergraphs. Later, Balogh and Samotij~\cite{BS20} developed an efficient container theorem, which they used to give a short proof of a result of Fox and Sudakov~\cite{FS09} on multicolour induced Ramsey numbers of graphs. More recently, a new container theorem of Campos and Samotij~\cite{CS24} was the main tool used by Aragão, Campos, Dahia, Filipe and Marciano~\cite{ACDFM25} to prove an exponential upper bound for diagonal induced Ramsey numbers. 

Applying this method (more precisely, a container theorem from~\cite{CS24}) in our setting, we will be able to find a small family of `containers' for pairs $G' \subset G$ of graphs with no copy of $F$ in $G'$ that is induced in $G$, each of which contains few cliques (see \Cref{thm:containers}). We will show, for each such container, that if there is also no induced copy of $H$ in the remaining edges of $G$, then $G$ avoids a fixed set of $\delta n^2$ edges of $K_n$. The theorem will then follow easily by bounding the probability of this event and using a union bound over the choice of the container. 

We remark that in order for this approach to work, we need to have at least (roughly) $n^{2-1/s}$ ``bits'' of randomness in our construction, in order to overcome the number of containers in the union bound. In~\cite{ACDFM25,BS20,CDFRS17}, like in most applications of the hypergraph container method to graphs, this is achieved by choosing the edges (or hyperedges) independently at random. In our setting we instead achieve it by taking the dimension of our affine space to be larger than $s$.

The rest of this paper is organised as follows. In \Cref{sec:construction} we describe our construction and state the container theorem (\Cref{thm:containers}) that we use to prove our main theorems. In \Cref{sec:supersaturation} we prove a simple supersaturation result (see \Cref{lem:dense-sets-cliques,lem:dense-sets-cliques:multicolour}), and in \Cref{sec:geometric} we exploit the geometric structure of our construction to find a set containing many blue edges (see \Cref{lem:blue-density}). In \Cref{sec:proof} we prove \Cref{thm:main,thm:multicolour}, and in \Cref{sec:remarks} we ask whether the  dependence on $v(F)$ of the exponent in~\eqref{eq:rFH:bound} can be improved from quadratic to linear. Finally, in \Cref{apx:containers}, we show how to deduce our container theorem from~\cite[Theorem D]{CS24}.

\medskip
\noindent \textbf{AI declaration and independent work:} The proofs of Theorems~\ref{thm:main} and~\ref{thm:multicolour} were developed by the authors without any use of AI tools, and the first draft of the proof was completed in May, 2026. The only use of AI tools in the writing was reviewing the final draft of the paper. On October 21st, 2026, an alternative proof of \Cref{thm:main}, using a different method, was posted to arXiv by  Zhong, Kashima and Mao~\cite{ZKM26}.

\section{Preliminaries}\label{sec:construction}

In this section we introduce the random geometric construction used in the proof and explain how we will use the container method to bound the probability that there exists a colouring that avoids given monochromatic induced subgraphs. Throughout the proof we will let $v(H) = k$ and $v(F_1) = \cdots = v(F_r) = s$, and assume that $s \ge 3$ and $e(H)>0$, since otherwise the desired bound is trivial. When $r = 1$ we may also assume that $k \ge s$, by swapping the roles of $F$ and $H$ if necessary, and noting that $x^{y^2} \le y^{x^2}$ for all $x \ge y \ge 3$. 

\subsection{The construction} 

Fix $r \in \N$ and let $q$ be a prime satisfying
\begin{equation}\label{def:q}
k^{6r^2s} \le q \le 2k^{6r^2s},
\end{equation}
and set
$$d = 4s \qquad \text{and} \qquad n = q^d.$$
Let $\cL$ be the set of affine lines in $\Fq^d$.  For each line $L\in\cL$, choose independently and uniformly at random a function
\begin{equation}\label{def:pi}
\pi_L : L \to V(H).
\end{equation}
For distinct points $x,y \in \Fq^d$, let $L_{xy}$ be the unique affine line containing $x$ and $y$. We now define a random graph $G = G(q,d,H)$ on vertex set $\Fq^d$ by setting
\begin{equation}\label{def:G}
xy \in E(G) \quad \Leftrightarrow \quad \pi_{L_{xy}}(x) \pi_{L_{xy}}(y) \in E(H).
\end{equation}
Thus, on each line $L$, the graph induced by $G$ is a random blow-up of $H$, where the parts are the fibres of $\pi_L$.

\subsection{The key container theorem}\label{sec:containers} 

The key tool in our proof is the following container theorem, which is very similar to the one proved by Balogh and Samotij~\cite{BS20}. In \Cref{apx:containers} we will show how to deduce it from a container theorem of Campos and Samotij~\cite[Theorem D]{CS24}.

\begin{theorem}\label{thm:containers}
Let $n \in \N$ and $\eps \in (0,1)$, and let $F$ be a graph with $s \ge 3$ vertices. 
There exists a family $\cC\subset 2^{E(K_n)} \times 2^{E(K_n)}$ of pairs of graphs satisfying
$$|\cC| \le \exp\Big( 8 s^4 \cdot \eps^{-\frac{2}{s+1}} \hspace{0.02cm} n^{2 - \frac{2}{s+1}} \log n \Big)$$
such that the following hold.
\begin{enumerate}[label=\textup{(\alph*)}, ref=(\alph*)]
\item \label{item:cont1} For every pair of graphs $G' \subset G$ with $n$ vertices for which there is no $U \subset V(G)$ satisfying $G'[U] = G[U] \cong F$, there exists $(X,Y) \in \cC$ such that
\begin{equation*}
G^c \subset X \qquad \text{and} \qquad G' \subset Y.
\end{equation*}
\item \label{item:cont2} For every $(X,Y) \in \cC$, there are at most $\eps n^s$ copies of $K_s$ in $X \cap Y$.
\end{enumerate}
\end{theorem}

Observe that if the graph $G$ admits a bad colouring, then the graph $G_B \subset G$ of blue edges does not contain any copies of $F$ which are induced in $G$. By \Cref{thm:containers}, it follows that there exists a container $(X,Y) \in \cC$ such that $G^c \subset X$, $G_B \subset Y$, and the graph $X \cap Y$ contains at most $\eps n^s$ copies of $K_s$. Our task will then be to show that 
$$e(X) < (1 - \delta) \binom{n}{2}$$
for some constant $\delta > 0$, which then easily implies that the event $G^c \subset X$ has sufficiently low probability that we may take a union bound over all the containers in $\cC$.

The key point in considering $X \cap Y$ is that this allows us to shift the perspective from induced copies of $F$ to cliques in $X \cap Y$. Then in case $ e(X)\geq (1-\delta)\binom{n}{2}$, we will be able to find many cliques in $X\cap Y$ by using the pseudorandomness of $G_B\subset Y$ and the fact $X$ is very dense. This will contradict \Cref{item:cont2} showing that for all containers we have $ e(X)< (1-\delta)\binom{n}{2}$.

\section{Supersaturation}\label{sec:supersaturation}

In this section we will prove a supersaturation lemma for $r$-colourings of pseudorandom graphs. The proof is simple, but to warm ourselves up, let's first prove it in the case $r = 1$, where the statement and proof are both a little easier to understand. Define
\begin{equation}\label{def:eps:mono}
\eps(s,\alpha) = \frac{\alpha^{\binom{s+1}{2}}}{s!}
\end{equation}
for each $s \in \N$ and $\alpha > 0$.

\begin{lemma}\label{lem:dense-sets-cliques}
Let $J$ be a graph on $n$ vertices, and let $s \in \N$ and $\alpha \in (0,1)$. If 
\begin{equation}\label{eq:dense-sets-cliques}
e\big( J[S] \big) \ge \alpha |S|^2
\end{equation}
for every $S \subset V(J)$ with $|S| \ge \alpha^{s-1} n$, then $J$ contains more than $\eps(s,\alpha) n^s$ copies of $K_s$.
\end{lemma}

\begin{proof}
Observe first that if~\eqref{eq:dense-sets-cliques} holds for a set $S$, then the graph $J[S]$ has more than $\alpha |S|$ vertices of degree at least $\alpha |S|$. We may therefore construct ordered $s$-cliques greedily, simply by choosing a vertex of degree at least $\alpha |S|$ in $J[S]$, where $S$ is the common neighbourhood of the vertices chosen so far. Since $|S| \ge \alpha^{s-1} n$ 
throughout the process, there are at least $\alpha |S|$ such vertices at each step. It follows that there are more than 
$$\prod_{i = 1}^{s} \alpha^i n = \alpha^{\binom{s+1}{2}} n^s$$
ordered $s$-cliques, and hence more than $\eps(s,\alpha) n^s$ copies of $K_s$, as claimed. 
\end{proof}

For the proof of \Cref{thm:multicolour}, we will instead need the following multicolour version of \Cref{lem:dense-sets-cliques}. Define
\begin{equation}\label{def:eps:multi}
\eps(s,r,\alpha) = \frac{\alpha^{\binom{t+1}{2}}}{t^s}
\end{equation}
for each $r,s \in \N$ and $\alpha > 0$, where $t = r(s-1) + 1$. 

\begin{lemma}\label{lem:dense-sets-cliques:multicolour}
Let $J$ be a graph on $n$ vertices, and let $r,s \in \N$ and $\alpha \in (0,1)$. If
\begin{equation}\label{eq:dense-sets-cliques:multicolour}
e\big( J[S] \big) \ge \alpha r |S|^2
\end{equation}
for every $S \subset V(J)$ with $|S| \ge \alpha^{r(s-1)} n$, then every $r$-colouring of the edges of $J$ contains more than $\eps(s,r,\alpha) n^s$ monochromatic copies of $K_s$.
\end{lemma}

\begin{proof}
Note that if~\eqref{eq:dense-sets-cliques:multicolour} holds for a set $S$, then there is a colour $i$ such that the graph $J[S]$ has more than $\alpha |S|$ vertices with degree at least $\alpha |S|$ in colour $i$. We may therefore construct ordered left-coloured cliques of size $r(s-1)+1$ greedily, simply by choosing such a vertex in each step, and intersecting $S$ with the colour $i$ neighbourhood of the chosen vertex. Since $|S| \ge \alpha^{r(s-1)} n$ throughout the process, there are at least $\alpha |S|$ such vertices at each step, and we therefore obtain at least
$$\prod_{i = 1}^{r(s-1)+1} \alpha^i n = \alpha^{\binom{t+1}{2}} n^t$$
ordered $t$-cliques, where $t = r(s-1)+1$. Each of these $t$-cliques contains at least one monochromatic copy of $K_s$, and each $K_s$ is included in at most $t^s n^{t-s}$ of the ordered $t$-cliques, so it follows that there are more than $\eps(s,r,\alpha) n^s$ monochromatic copies of $K_s$, as claimed. 
\end{proof}

\section{Geometric lemmas}\label{sec:geometric}

In this section we will prove the following lemma.

\begin{lemma}\label{lem:blue-density}
With high probability, in every red-blue colouring of the edges of $G$ with no red induced copy of $H$, we have
$$e(G_B[S])\ge \frac{|S|^2}{2^6 k^2}$$
for every set $S \subset V(G)$ of size $|S| \ge 2k^2 q^{d-1}$. 
\end{lemma}

To prove this, we will need two lemmas that exploit the geometric structure of our construction. The first of these says that every sufficiently large subset of $V(G) = \F_q^d$ has a large intersection with at least half of the lines in $\cL$.

\begin{lemma}\label{lem:many-lines} 
For every set $S \subset V(G)$ of size $|S| \ge 8 q^{d-1}$, we have
$$|S \cap L| \ge \frac{|S|}{2q^{d-1}}$$
for at least half of the lines $L \in \cL$.
\end{lemma}

\begin{proof}
Let $L \in \cL$ be a uniformly chosen random line, and observe that 
$$\Ex\big[ | S \cap L | \big] = \frac{|S|}{q^{d-1}}$$
and
$$\Var\big( |S \cap L| \big) = \sum_{x,y \in S} \Pr\big( x,y \in L \big) - q^{-2(d-1)} |S|^2.$$
Since $\Pr(x \in L ) = q^{-(d-1)}$ and $\Pr\big( x,y \in L \big) = \frac{1}{|\cL|} \le q^{-2(d-1)}$ for every $x \ne y$, it follows that
$$\Var\big( |S \cap L| \big) \le \frac{|S|}{q^{d-1}}.$$
Therefore, by Chebyshev's inequality, we obtain
\begin{equation*}
\Pr\bigg( | S \cap L | < \frac{\Ex\big[ |S \cap L| \big]}{2} \bigg) \le 4 \cdot \frac{\Var\big( |S \cap L| \big)}{\Ex\big[ |S \cap L |\big]^2} \le \frac{4q^{d-1}}{|S|} \le \frac{1}{2},
\end{equation*}
as claimed, since $|S| \ge 8 q^{d-1}$.
\end{proof}

\Cref{lem:many-lines} guarantees that every large subset of $V(G)$ has substantial intersections
with many lines. Together with the random choice of the maps $\pi_L$, this will imply
that large subsets span many edges of $G$. This allows us to show that, with high probability, for many of these lines $L\in\cL$ satisfying \Cref{lem:many-lines}, the corresponding blow-up of $H$ is balanced even after restriction to a large set. To be precise, we have the following lemma.

\begin{lemma}\label{lem:balanced-lines}
With high probability, for every set $S \subset V(G)$ of size $|S| \ge 2k^2 q^{d-1}$, at least $|\cL|/4$ lines $L \in \cL$ satisfy
\begin{equation}\label{eq:BalancedPartition}
\big| \pi_L^{-1}(u) \cap S \big| \ge \frac{|S|}{4kq^{d-1}}
\end{equation}
for every $u \in V(H)$.
\end{lemma}

\begin{proof}
Fix a set $S \subset V(G)$ of size $|S| \ge 2k^2 q^{d-1}$, and let $\cL_S$ be the family of lines given by \Cref{lem:many-lines} applied to $S$, so $|\cL_S| \ge |\cL| / 2$ and
$$|S \cap L| \ge \frac{|S|}{2q^{d-1}} \ge k^2$$ 
for every $L \in \cL_S$. Fix $L \in \cL_S$ and $u \in V(H)$, and observe that the random variable 
$$Z = \big| \pi_L^{-1}(u) \cap S \big|$$ 
has a binomial distribution with mean $|S \cap L| / k$.  By Chernoff's inequality, it follows that
$$\Pr\bigg( \big|\pi_L^{-1}(u) \cap S \big| < \frac{|S|}{4kq^{d-1}} \bigg) \le \Pr\bigg( \big| \pi_L^{-1}(u) \cap S \big| < \frac{|S \cap L|}{2k} \bigg) \le e^{-k/8}.$$
Taking a union bound over $u \in V(H)$, we obtain
$$\Pr\bigg( \min_{u \in V(H)} \big| \pi_L^{-1}(u) \cap S \big| < \frac{|S|}{4kq^{d-1}} \bigg) \le k \cdot e^{-k/8} < \frac{1}{4},$$
since $k$ is sufficiently large. 

Now, recall that the random functions $\pi_L$ on distinct lines $L \in \cL$ are independent. Therefore, by another Chernoff bound, we may bound the probability that~\eqref{eq:BalancedPartition} fails to hold for more than half of the lines in $\cL_S$ by
$$\Pr\bigg( \bigg| \bigg\{ L \in \cL_S :~\min_{u \in V(H)} \big| \pi_L^{-1}(u) \cap S \big| < \frac{|S|}{4kq^{d-1}} \bigg\} \bigg| \ge \frac{|\cL_S|}{2} \bigg) \le \exp\bigg( - \frac{|\cL_S|}{8} \bigg).$$
Since $|\cL_S| \ge \frac{1}{2} \cdot q^{2d-2}$, a union bound over all choices of $S$ gives a bound on the probability that such a set $S$ exists of 
$$2^n \exp\bigg( - \frac{q^{2d-2}}{16} \bigg) \le \exp\bigg( q^d - \frac{q^{2d-2}}{16} \bigg) = o(1),$$
as claimed.
\end{proof}

We next recall the following key observation from~\cite{KPR98}. For the reader's convenience, we give the short proof. Recall that we write $H(t)$ to denote the $t$-blow-up of $H$. 

\begin{observation}[Lemma~13 of~\cite{KPR98}]\label{obs:ksquared}
A red-blue colouring of $H(t)$ with no red induced copy of $H$ contains at least $t^2$ blue edges.
\end{observation} 

\begin{proof}
There are $t^k$ transversal copies of $H$ in $H(t)$, and each must have at least one blue edge. Since each edge is contained in $t^{k-2}$ transversal copies of $H$, the claim follows.
\end{proof}

Combining this observation with \Cref{lem:balanced-lines}, we obtain the claimed lower bound on the density of blue edges in a large set of vertices. 

\begin{proof}[Proof of \Cref{lem:blue-density}]
Let $\cA$ denote the event that the conclusion of \Cref{lem:balanced-lines} holds: that for every set $S \subset V(G)$ of size $|S| \ge 2k^2 q^{d-1}$, at least $|\cL|/4$ lines $L \in \cL$ satisfy
\begin{equation}\label{eq:BalancedPartition:again}
\big| \pi_L^{-1}(u) \cap S \big| \ge t = \frac{|S|}{4kq^{d-1}}
\end{equation}
for every $u \in V(H)$. By \Cref{lem:balanced-lines}, we may assume that $\cA$ holds. 

Fix a red-blue colouring of the edges of $G$ with no red induced copy of $H$, and a set $S \subset V(G)$ of size $|S| \ge 2k^2 q^{d-1}$. We claim that the desired bound 
$$e(G_B[S])\ge \frac{|S|^2}{2^6 k^2}$$
follows deterministically from the event $\cA$. To prove this, let $\cL_S \subset \cL$ be the set of lines given by the event $\cA$, so $|\cL_S| \ge |\cL|/4$ and~\eqref{eq:BalancedPartition:again} holds for every $L \in \cL_S$ and every $u \in V(H)$. For each $L \in \cL_S$ and $u \in V(H)$, choose a subset
$$T_u \subset \pi_L^{-1}(u) \cap S \qquad \text{with} \qquad |T_u| = t,$$
and set 
$$T = \bigcup_{u \in V(H)} T_u.$$ 
The construction of $G$ implies that $G[T] = H(t)$: every pair with one endpoint in $T_u$ and the other in $T_v$ is an edge of $G$ exactly when $uv \in E(H)$.

By \Cref{obs:ksquared}, it follows that $G[T] \subset G[S]$ has at least $t^2$ blue edges. Since every pair of points lies on a unique line and $|\cL_S| \ge |\cL|/4 \ge q^{2d-2}/4$, we conclude that
$$e(G_B[S]) \ge t^2 \cdot |\cL_S| \ge \bigg( \frac{|S|}{4kq^{d-1}} \bigg)^2 \cdot \frac{q^{2d-2}}{4} = \frac{|S|^2}{2^6 k^2},$$
as claimed.
\end{proof}

\section{Proof of the main theorems}\label{sec:proof}

In this section, we will prove \Cref{thm:main,thm:multicolour}. To do so, we will first apply \Cref{thm:containers} to find for each colouring of $G$ with no blue induced copy of $F$ and no red induced copy of $H$ a container pair $(X,Y)$.

We will then use \Cref{lem:dense-sets-cliques} to show that there exists a large subset $S$ where $X \cap Y$ is sparse. At the same time \Cref{lem:blue-density} will show that $G_B$, and therefore $Y\supset G_B$, is dense in $S$, which will imply an upper bound on $e(X)$. The final argument is then to union bound over pairs $(X,Y)$ and use the property that $G^c \subset X$.

\begin{proof}[Proof of \Cref{thm:main}] 
Let $F$ and $H$ be graphs with $v(F) = s$ and $v(H) = k$. We may assume that $k \ge s \ge 3$ and $k \ge 2^7$, and also that $e(H) > 0$, since otherwise the claimed bound holds trivially. Set $\alpha = k^{-3}$, define $\eps = \eps(s,\alpha)$ as in~\eqref{def:eps:mono}, and apply \Cref{thm:containers} to the graph $F$ with $n = q^d$ and $\eps = \eps(s,\alpha)$. We obtain a family $\cC$ of pairs of graphs such that 
\begin{equation}\label{eq:size:of:C}
|\cC| \le \exp\Big( 8s^4 \hspace{0.01cm} \eps^{-\frac{2}{s+1}} \hspace{0.02cm} n^{2 - \frac{2}{s+1}} \log n \Big)
\end{equation}
and properties \ref{item:cont1} and \ref{item:cont2} in \Cref{thm:containers} hold.   

Let $G$ be the random graph constructed in \Cref{sec:construction}, with $r = 1$, and suppose that there exists a red-blue colouring of $G$ with no blue induced copy of $F$ and no red induced copy of $H$. Recall from \Cref{lem:blue-density} that with high probability we have
\begin{equation}\label{eq:eBS:lower:bound}
e(G_B[S]) \ge \frac{|S|^2}{2^6 k^2} \ge 2\alpha |S|^2
\end{equation}
for every set $S \subset V(G)$ of size $|S| \ge 2k^2 q^{d-1}$ in every red-blue colouring of the edges of $G$ with no red induced copy of $H$. Since the graph of blue edges $G_B \subset G$ is induced $F$-free in $G$, property~\ref{item:cont1} of \Cref{thm:containers} implies that there exists $(X,Y) \in \cC$ such that
$$G^c \subset X \qquad \text{and} \qquad G_B \subset Y.$$
Set $\delta = \alpha^{2s-1}$. We make the following claim. 

\begin{claim}\label{claim:eX:small}
$$e(X) \le \big( 1 - \delta \big)\binom{n}{2}.$$
\end{claim}

\begin{proof}[Proof of \Cref{claim:eX:small}]
\renewcommand{\qedsymbol}{$\blacksquare$}
By \Cref{lem:dense-sets-cliques} and property~\ref{item:cont2} of the container family $\cC$, there exists a set $S \subset V(G)$ with $|S| \ge \alpha^{s-1} n$
that contains at most $\alpha |S|^2$ edges of the graph $X \cap Y$. Noting that
\begin{equation}\label{eq:S:lower:bound}
|S| \ge \alpha^{s-1} n = \frac{q^d}{k^{3s-3}} \ge 2k^2 q^{d-1},
\end{equation}
since $q \ge k^{3s}$, it follows from~\eqref{eq:eBS:lower:bound} that
$$e(G_B[S]) \ge 2\alpha |S|^2.$$
Since $G_B \subset Y$, at least $\alpha |S|^2$ of these edges are not contained in $X$, and hence 
$$e(X) \le {n \choose 2} - \alpha |S|^2 \le \big( 1 - \delta \big) \binom{n}{2},$$
by~\eqref{eq:S:lower:bound} and our choice of $\delta = \alpha^{2s-1}$.
\end{proof}

It now only remains to bound the probability of the event $G^c \subset X$, and take a union bound over the choice of the container $(X,Y)$. To do so, we choose a maximal collection $W \subset X^c$ of non-edges of $X$ that are contained in pairwise distinct lines of $\Fq^d$. Since each line contains at most $\binom{q}{2}$ non-edges of $X$, it follows from \Cref{claim:eX:small} that 
$$|W| \ge \delta {n \choose 2} {q \choose 2}^{-1} \ge \frac{\delta n^2}{q^2}.$$
Now, recall from~\eqref{def:G} that for each pair $xy \in W$, we have
$$xy \in E(G) \quad \Leftrightarrow \quad \pi_{L_{xy}}(x) \pi_{L_{xy}}(y) \in E(H),$$
and therefore, by the definition~\eqref{def:pi} of $\pi_{L_{xy}}$, 
$$\Pr\big( xy \in E(G) \big) = \frac{2e(H)}{k^2} \le 1 - \frac{1}{k}.$$
Since the maps $\pi_L$ are chosen independently for distinct lines, the events $\{ xy \in E(G) \}$ are independent whenever the pairs $xy$ lie on distinct lines, so 
$$\Pr\big( G^c \subset X \big) \le \Pr\big( W \subset G \big) \le \bigg( 1 - \frac{1}{k} \bigg)^{\delta n^2 / q^2} \le \exp\bigg( - \frac{\alpha^{2s} n^2}{q^2} \bigg).$$
Finally, taking a union bound over the family $\cC$ of containers, and recalling that~\eqref{eq:eBS:lower:bound} holds with high probability by \Cref{lem:blue-density}, we deduce that
\begin{equation}\label{eq:bound}
\Pr\big( G \nindarrow (F,H) \big) \le |\cC| \cdot\exp\bigg( - \frac{\alpha^{2s} n^2}{q^2} \bigg) + o(1).
\end{equation}

Now, since $\eps = \alpha^{\binom{s+1}{2}}/s!$ and $n = q^{4s}$, we have
$$\eps^{-\frac{2}{s+1}} \le \frac{s^2}{\alpha^s} \qquad \text{and} \qquad n^{-\frac{2}{s+1}} \log n \le \frac{4s}{q^5},$$
since $s! \le s^{s+1}$, $\log q \le q$ and $s \ge 3$, and hence, by~\eqref{eq:size:of:C},  
\begin{equation*}
\frac{\log |\cC|}{n^2} \le 8s^4 \eps^{-\frac{2}{s+1}} n^{- \frac{2}{s+1}} \log n \le \frac{2^5 s^7}{\alpha^s q^5}.
\end{equation*}
Moreover, since $q \ge k^{4s}$ and $\alpha = k^{-3}$, we have $q^3 \alpha^{3s} \ge k^{3s} \ge 2^6 s^7$, and hence 
\begin{equation*}
\frac{\log |\cC|}{n^2} \le \frac{2^5 s^7}{\alpha^s q^5} \le \frac{\alpha^{2s}}{2q^2},
\end{equation*}
Together with~\eqref{eq:bound}, this implies that
$$\Pr\big(G \nindarrow (F,H)\big) \le \exp\bigg(-\frac{\alpha^{2s}n^2}{2q^2}\bigg) + o(1) = o(1).$$
Thus $G \indarrow (F,H)$ with high probability, and hence, by~\eqref{def:q}, we have
$$r_{\mathrm{ind}}(F,H) \le n = q^{4s} \le 2^{4s} k^{24s^2} \le k^{25s^2},$$
as required.
\end{proof}

The proof of \Cref{thm:multicolour} is almost identical: we simply replace the application of \Cref{lem:dense-sets-cliques} with an application of \Cref{lem:dense-sets-cliques:multicolour}. 

\begin{proof}[Proof of \Cref{thm:multicolour}] 
Let $F_1,\ldots,F_r$ and $H$ be graphs with $v(F_1) = \cdots = v(F_r) = s$ and $v(H) = k$. We may  assume that $s \ge 3$, that $k \ge r + 2^7$ (by our assumption that $v(H) \ge Cr$) and that $e(H) > 0$. Set $\alpha = r^{-1} k^{-3}$, define $\eps = r^{-1} \eps(s,r,\alpha)$, where $\eps(s,r,\alpha)$ was defined in~\eqref{def:eps:multi}, and apply \Cref{thm:containers} to $F_i$ for each $i \in [r]$. We obtain families $\cC_1,\ldots,\cC_r$ of containers, with
\begin{equation}\label{eq:size:of:Ci}
    |\cC_i| \le \exp\Big( 8s^4 \hspace{0.02cm} \eps^{-\frac{2}{s+1}} \hspace{0.02cm} n^{2 - \frac{2}{s+1}} \log n \Big),
\end{equation}
such that properties \ref{item:cont1} and \ref{item:cont2} in \Cref{thm:containers} hold for each $\cC_i$ (with respect to $F_i$). 

Let $G$ be the random graph constructed in \Cref{sec:construction}, and recall from \Cref{lem:blue-density} that with high probability we have
\begin{equation}\label{eq:eBS:lower:bound:r}
e(G_B[S]) \ge \frac{|S|^2}{2^6 k^2} \ge 2\alpha r |S|^2
\end{equation}
for every set $S \subset V(G)$ of size $|S| \ge 2k^2 q^{d-1}$ and every red-blue colouring of the edges of $G$ with no red induced copy of $H$. Suppose that there exists an $(r+1)$-colouring of $G$ with no induced copy of $F_i$ in colour $i$ for each $i \in [r]$ and no induced copy of $H$ in colour $r+1$. By property~\ref{item:cont1} of \Cref{thm:containers}, for each $i \in [r]$ there exists a container $(X_i,Y_i) \in \cC_i$ such that
$$G^c \subset X_i \qquad \text{and} \qquad G_i \subset Y_i,$$
where $G_i$ denotes the set of edges of colour $i$. We make the following claim. 

\begin{claim}\label{claim:eX:small:multi}
Set $\delta = \alpha^{2rs-1}$. Then
$$e(X_i) \le \big( 1 - \delta \big)\binom{n}{2}$$
for some $i \in [r]$. 
\end{claim}

\begin{proof}[Proof of \Cref{claim:eX:small:multi}]
\renewcommand{\qedsymbol}{$\blacksquare$}
Define $J$ to be the graph with vertex set $V(G)$ and edge set 
$$E(J) = \bigcup_{i = 1}^r E\big( X_i \cap Y_i \big).$$ 
By \Cref{lem:dense-sets-cliques:multicolour} and property~\ref{item:cont2} of the container families $\cC_1,\ldots,\cC_r$, there exists $S \subset V(G)$ with $|S| \ge \alpha^{r(s-1)} n$ such that 
$$e\big( J[S] \big) < \alpha r |S|^2.$$
Noting that
\begin{equation}\label{eq:S:lower:bound:multi}
|S| \ge \alpha^{r(s-1)} n = \frac{q^d}{(r k^3)^{r(s-1)}} \ge 2k^2 q^{d-1},
\end{equation}
since $q \ge k^{4r^2s} \ge (r k^3)^{rs}$, it follows from \Cref{lem:blue-density} that
$$e(G_B[S]) \ge 2\alpha r |S|^2,$$
where blue is the union of the colours in $[r]$. Since $G_B \subset Y_1 \cup \cdots \cup Y_r$, it follows that there exists $i \in [r]$ such that at least $\alpha |S|^2$ of the edges of $Y_i$ are not contained in $X_i$, and hence 
$$e(X_i) \le {n \choose 2} - \alpha |S|^2 \le \big( 1 - \delta \big) \binom{n}{2},$$
by~\eqref{eq:S:lower:bound:multi} and our choice of $\delta = \alpha^{2rs-1}$.  
\end{proof}

It now only remains to bound the probability of the event $G^c \subset X_i$, and take a union bound over the choice of the container $(X_i,Y_i)$. To do so, we choose a collection $W \subset X_i^c$ of $\delta n^2 / q^2$ non-edges of $X_i$ that are contained in pairwise distinct lines of $\Fq^d$, and calculate as in the proof of \Cref{thm:main}, to obtain
$$\Pr\big( G^c \subset X_i \big) \le \Pr\big( W \subset G \big) \le \bigg( 1 - \frac{1}{k} \bigg)^{\delta n^2 / q^2} \le \exp\bigg( - \frac{\alpha^{2rs} n^2}{q^2} \bigg).$$
Taking a union bound over the families $\cC_1, \ldots, \cC_r$ of containers, and recalling that~\eqref{eq:eBS:lower:bound:r} holds with high probability by \Cref{lem:blue-density}, we deduce that
\begin{equation}\label{eq:bound:r}
\Pr\big( G \nindarrow (F_1,\ldots,F_r,H) \big) \le \sum_{i = 1}^r |\cC_i| \cdot \exp\bigg( - \frac{\alpha^{2rs} n^2}{q^2} \bigg) + o(1).
\end{equation}
To bound $|\cC_i|$, recall that $\eps = \alpha^{\binom{t+1}{2}}/(rt^s)$, and note that $\binom{t+1}{2} \le r^2s(s+1)/2$, since $t \le rs$. It follows that
$$\eps^{-\frac{2}{s+1}} \le \frac{(rs)^2}{\alpha^{r^2s}} \qquad \text{and} \qquad n^{-\frac{2}{s+1}} \log n \le \frac{4s}{q^5},$$
since $n = q^{4s}$ and $s \ge 3$, and therefore, by~\eqref{eq:size:of:Ci},  
\begin{equation*}
\frac{\log |\cC_i|}{n^2} \le 8s^4 \hspace{0.02cm} \eps^{-\frac{2}{s+1}} \hspace{0.02cm} n^{- \frac{2}{s+1}} \log n \le \frac{1}{\alpha^{2r^2s} q^5},
\end{equation*}
since $\alpha = r^{-1} k^{-3}$. Noting that $q^3 \ge k^{18r^2s} \ge 2\alpha^{-4r^2s}$, since $k \ge r$, it follows that
\begin{equation*}
\frac{\log |\cC_i|}{n^2} \le \frac{1}{\alpha^{2r^2s} q^5} \le \frac{\alpha^{2r^2s}}{2q^2},
\end{equation*}
Together with~\eqref{eq:bound:r}, this implies that
$$\Pr\big( G \nindarrow (F_1,\ldots,F_r,H) \big) \le r \exp\bigg( - \frac{\alpha^{2rs} n^2}{2q^2} \bigg)+o(1) = o(1).$$
Thus $G \indarrow (F_1,\ldots,F_r,H)$ with high probability, and hence, by~\eqref{def:q}, we have  
$$r_{\mathrm{ind}}(F_1,\ldots,F_r,H) \le n = q^{4s} \le 2^{4s} k^{24r^2s^2} \le k^{25r^2s^2},$$
as required.
\end{proof}

\section{Concluding remarks}\label{sec:remarks}

We have proved that, for every fixed graph $F$, the off-diagonal induced Ramsey number $r_{\mathrm{ind}}(F,H)$ is polynomial in $v(H)$.  This confirms the conjecture of Kohayakawa, Pr\"omel, and R\"odl in its qualitative form.  Quantitatively, however, the dependence of the exponent on $v(F)$ remains to be optimized: our argument gives
\[r_{\mathrm{ind}}(F,H) \leq v(H)^{O(v(F)^2)}.\]
In particular whenever $v(F)\geq \log v(H)$, the bound proved by Kohayakawa, Pr\"omel, and R\"odl is better than ours since they prove \[r_{\mathrm{ind}}(F,H) \leq v(H)^{O(v(F) \log v(H))}.\]

A natural next problem is to then determine the correct dependence of the exponent on the fixed graph $F$. This leads to the following question.

\begin{question}
Is it true that 
$$r_{\mathrm{ind}}(F,H) \le v(H)^{O(v(F))}$$
for all graphs $F$ and $H$?
\end{question}

If true, this would give the optimal form of the polynomial dependence on $v(H)$. In general the best lower bound on $r_{\mathrm{ind}}(F,H)$ that we are aware of is just the usual Ramsey number $r(F,H)$. It would also be interesting to improve this lower bound.

\section*{Acknowledgements}

The authors are grateful to Rob Morris for rewriting the article to improve the presentation and for pointing out how our proof could be generalized to more than two colours.

\bibliographystyle{abbrv}
\def\bibfont{\footnotesize}
\bibliography{refs}

\appendix
\crefalias{section}{appendix}
\crefname{appendix}{appendix}{appendices}

\section{Proof of \Cref{thm:containers}}\label{apx:containers}

We will use the following container theorem of Campos and Samotij~\cite{CS24} to deduce \Cref{thm:containers}.

\begin{theorem}[\hspace{-0.2mm}\cite{CS24}, Theorem D]\label{thm:CSTheoremD}
Let $\mathcal{H}$ be an $r$-uniform hypergraph with a finite vertex set $V$. For every $p \in (0, 1/(8r^2)]$, there exists a family $\mathcal{S} \subset 2^V$ and functions
\[
g:\mathcal{I}(\mathcal{H}) \to \mathcal{S}
\qquad \text{and} \qquad
f:\mathcal{S} \to 2^V
\]
such that:
\begin{enumerate}[label=(\alph*)]
    \item For each $I \in \mathcal{I}(\mathcal{H})$, we have $g(I) \subset I \subset f(g(I))$.\label{item:1}

    \item Each $S \in \mathcal{S}$ has at most $8r^2p|V|$ elements.\label{item:2}

    \item For every $S \in \mathcal{S}$, there does not exist a hypergraph $\mathcal{H}_S \subset \mathcal{H}[f(S)]$ that satisfies
    \[
    \Delta_\ell(\mathcal{H}_S)
    <
    \frac{p^{\ell-1}e(\mathcal{H}_S)}{|f(S)|}
    \qquad
    \text{for all } \ell \in \{2,\dots,r\}.
    \]\label{item:3}
\end{enumerate}
\end{theorem}

We will apply \Cref{thm:CSTheoremD} to the following hypergraph, which encodes induced copies of $F$ in a subgraph of $K_n$. Let $\cH$ be a hypergraph with vertex set
\[V(\cH)=E(K_n)\times\{0,1\},\]
and whose edges are the sets of the form
\[
    \big(\varphi(E(F))\times\{1\}\big)\cup\big(\varphi\big(E(K_s)\setminus E(F)\big)\times\{0\}\big),
\]
where $\varphi: V(F) \to V(K_n)$ is an arbitrary embedding. In other words, $\cH$ encodes all the possible induced copies of $F$ in a subgraph of $K_n$. In particular, every pair of graphs $G_0 \subset G_1$ on $n$ vertices for which there is no $U\subset V(G_0)$ satisfying $G_0[U]=G_1[U]\cong F$ can be viewed as an independent set in $\cH$ considering the set
\[I=\big\{(e,1):e\in E(G_0)\big\} \cup \big\{(e,0):e\in E(\overline{G_1})\big\} \subset V(\cH).\]

With these observations in mind, we are now ready to prove \Cref{thm:containers}.

\begin{proof}[Proof of \Cref{thm:containers}] 
Let $p = (\eps n)^{-\frac{2}{s+1}}$ and $r = \binom{s}{2}$. For a set $A\subset V(\cH)$ and $i\in\{0,1\}$, define
\[A^{(i)} = \big\{ e \in E(K_n) : (e,i)\in A \big\}.\] 

First we deal with the trivial case that $p > 1/(32r^2)$. In this case we can take 
\[\cC = \big\{ (I^{(0)}, I^{(1)}): I \in \cI(\cH) \big\},\]
the family of all independent sets of $\cH$. Note that property~\ref{item:cont1} of \Cref{thm:containers} holds immediately, and for each $I \in \cI(\cH)$ there is no copy of $K_s$ in $I^{(0)} \cap I^{(1)}$, since $I$ spans no edge, implying property~\ref{item:cont2} in \Cref{thm:containers}. Moreover, we have
\[
    |\mathcal I(\cH)|\leq 2^{|V(\cH)|}
    \leq 2^{n^2}
    \leq
    \exp\big(8s^4
    \eps^{-\frac{2}{s+1}}n^{2-\frac{2}{s+1}}\log n\big),
\]
since $8s^4\eps^{-\frac{2}{s+1}} n^{-\frac{2}{s+1}} \geq 1$.

We may now assume $p \le 1/(32r^2)$. Applying \Cref{thm:CSTheoremD} to the $r$-uniform hypergraph $\cH$, we obtain a family $\cS \subset 2^{V(\cH)}$ and functions
\[
g:\mathcal{I}(\mathcal{H}) \to \mathcal{S}
\qquad \text{and} \qquad
f:\mathcal{S} \to 2^{V(\cH)}
\]
satisfying \cref{item:1,item:2,item:3} in \Cref{thm:CSTheoremD}.

Then define
\[\cC=\big\{(C^{(0)},C^{(1)}):~C\in f(\cS)\big\}.\]
Fix a pair of graphs $G_0 \subset G_1$ on $n$ vertices for which there is no $U\subset V(G_1)$ satisfying $G_0[U]=G_1[U]\cong F$. Note that the set
\[I=\big\{(e,1):e\in E(G_0)\big\} \cup \big\{(e,0):e\in E(G_1^c)\big\} \subset V(\cH)\]
is an independent set in $\cH$.
By \cref{item:1} in \Cref{thm:CSTheoremD}, there exists $S \in \cS$ such that $S=g(I)$ and $S\subset I \subset f(S)=C$. In particular, we must have
\[I^{(1)}=G_0\subset C^{(1)}\qquad\text{and}\qquad I^{(0)}=G_1^c\subset C^{(0)},\]
which proves property~\ref{item:cont1} of \Cref{thm:containers}. Moreover, by \cref{item:2} in \Cref{thm:CSTheoremD} the size of the family $\cC$ can be bounded as

\[|\cC|\le \binom{n^2}{\leq 8r^2pn^2} \le \exp \left( 8 s^4\cdot \eps^{-\frac{2}{s+1}} n^{2-\frac{2}{s+1}}\log n\right)\, ,\]
since $|V|\leq n^2 $ and using that $p = (\eps n)^{-\frac{2}{s+1}}$.

It remains to prove \cref{item:cont2} in \Cref{thm:containers}. Given $(X,Y)\in \cC$ let 
$$\cK_{X,Y}:=\left\{K\subset X\cap Y\, :\, K \cong K_s\right\}\, , $$ 
that is a hypergraph encoding copies of $K_s$ in $X \cap Y$. We have that $e(\cK_{X,Y})$ is the number of copies of $K_s$ in $X\cap Y$, therefore we only need to show that $e(\cK_{X,Y})\leq \eps n^s$.

Let $C=(X\times \{0\}) \cup (Y\times \{1\})$ and take some injective function $\phi: \cK_{X,Y} \to \cH[C]$ such that $\phi(K)\subset K\times \{0,1\}$, for each $K\in \cK_{X,Y}$. Let $\cH_C:= \phi(\cK_{X,Y})$, by \cref{item:3} in \Cref{thm:CSTheoremD}, this hypergraph must satisfy
$$e(\cH_C) \le \dfrac{|C|\Delta_\ell(\cH_C)}{p^{\ell-1}}$$
for some $\ell \in \big\{2,3,\ldots,\binom{s}{2}\big\}$.
First, we bound $\Delta_\ell(\cH_C)$. Let $v_\ell\in\{3,\ldots,s\}$ be the unique integer satisfying $\binom{v_\ell - 1}{2}<\ell\le\binom{v_\ell}{2}$. Thus,
\begin{equation}\label{eq:2}
\Delta_{\ell}(\cH_C)\le \Delta_\ell(\cK_{X,Y})\le \binom{n-v_\ell}{s-v_\ell}\leq n^{s-v_\ell},
\end{equation}
where in the first inequality we use that for each $E\in \cH_C$ there is some $K\in \cK_{X,Y}$ such that $E\subset K$, and for the second we use that any $\ell$ pairs span at least $v_\ell$ vertices. 

Now observe that
\[\ell-1\le\binom{v_\ell}{2}-1=\frac{(v_\ell+1)(v_\ell-2)}{2}\le\frac{(s+1)(v_\ell-2)}{2},\]
which implies
\begin{equation*}
    \frac{2(\ell-1)}{s+1}\le v_\ell-2 \, .
\end{equation*}
Then we obtain 
\begin{equation} \label{eq:3}
p^{-\ell+1}=(\eps n)^{\frac{2(\ell-1)}{s+1}}\leq (\eps n)^{\frac{2(\ell-1)}{v_\ell+1}}\leq (\eps n)^{v_\ell-2}.
\end{equation}

Putting together \eqref{eq:2} and \eqref{eq:3} we obtain $$\frac{|C|\Delta_\ell(\cH_C)}{p^{\ell-1}}\leq n^{s-v_\ell +2}(\eps n)^{v_\ell-2} \leq \eps n^s \, , $$ using that $v_\ell\geq 3$ and $\eps \in (0,1)$.
This finishes the proof since $$e(\cK_{X,Y})=e(\cH_C) \le \eps n^s.$$
\end{proof}
\end{document}